\documentclass[12pt]{article}

\newtheorem{lem}{Lemma}[section]
\newtheorem{prop}{Proposition}[section]
\newtheorem{thm}{Theorem}[section]
\newtheorem{cor}{Corollary}[section]

\newenvironment{pf}
{\noindent {\sc Proof.}} {$\hfill\rlap{$\sqcup$}\sqcap$\bigskip}

\let\cj=\overline
\def\f#1{{\mathcal{F}}_{#1}}
\def\pt#1{\left(#1\right)}

\begin{document}

\title{Bounds on the degree of APN polynomials\\ The Case of $x^{-1}+g(x)$}
\author{Gregor Leander\thanks{Faculty of Mathematics, Technical University of Denmark, {e-mail:} {\tt g.leander@mat.dtu.dk}}
\and
Fran\c cois Rodier\thanks{Institut of Mathematiques of Luminy --
C.N.R.S. --
{e-mail:} {\tt rodier@iml.univ-mrs.fr}}}
\date{}
\maketitle

\begin{abstract}
We prove that functions $f:\f{2^m} \to \f{2^m}$ of the form $f(x)=x^{-1}+g(x)$ where $g$ is any non-affine polynomial are APN on at most a finite number of fields $\f{2^m}$. Furthermore we prove that when the degree of $g$ is less then $7$ such functions are APN only if $m \le 3$ where these functions are equivalent to $x^3$.
\end{abstract}
\section{Introduction}
For a given integer $m$ denote by $q=2^m$ and let $\f{q}$ be the finite field with $q$ elements. We study functions
$f:\f{q}\to \f{q}$ given by their polynomial representation. Such a function, or polynomial, is called {\em{almost perfect nonlinear}} (APN)
if for every non-zero $a\in \f{q}$ and every  $b\in \f{q}$ the equation
\[ f(x)+f(x+a)=b\]
admits at most two solutions $x\in \f{2^m}$. Amongst others, APN functions have applications in cryptography. Namely, when used as an S-box in a block cipher being APN
ensures a good resistance against differential attacks.

Until 2006, there where only very few APN functions known and all of them where power mappings. It was even conjectured that any APN function is equivalent to one of the known APN power functions.
Here equivalence is usually defined by saying that two functions $f,g:\f{2^m}\to \f{2^m}$ are equivalent if there exist an affine permutation on $\f{2^m} \times \f{2^m}$ such that
the graph of $f$, i.e. the set $\{(x,f(x)) \ : \ x \in \f{2^{m}}\}$, is mapped to the graph of $g$. This equivalence is called CCZ-equivalence (see \cite{DBLP:journals/dcc/CarletCZ98})
and preserves the APN property.

In \cite{DBLP:journals/IEEE/Pot06} the first APN function which was not equivalent to any power function was found. Shortly after this, several infinite families of APN functions have been discovered, see for example \cite{bracken2008fmq,bracken2007nfq,budaghyan2008cqa,TraceAPN,budaghyan2008tcq}.

The problem of classifying all APN functions seems elusive today. Even the problem of classifying all APN power functions is an open problem and not a lot progress has been made here.
However, there are possible steps that can be taken towards a full classification. One approach that already proved to be successful is to show that certain polynomials
are not APN for infinitely many extensions of $\f{2}$. So here one first fixes a finite field $\f{q}$ and a function $f:\f{q}\to \f{q}$ given as a polynomial in $\f{q}[x]$ and
poses the question if this function can be APN on infinitely many extensions of $\f{q}$. There is a variety of classes of functions for which it can be shown
that each function is APN at most for a finite number of extensions. For example, Jedlicka \cite{Jedlicka} studied the case of power functions and Voloch \cite{Voloch} focused on binomials.
Moreover, Rodier studied general polynomials with the same approach in \cite{rodier2006bdp,rodierBFCA08}.

\section{Our Results}
So far the above question has only be treated when the polynomial description of $f$ does not depend on $n$. In this paper we tackle the question if a given class of functions can be APN infinitely many often
for functions of the form
\begin{eqnarray*}
f: \f{q} &\to& \f{q} \\
 f(x)&=&x^{q-2}+g(x)
\end{eqnarray*}
where $g$ is a fixed polynomial. The description of those functions clearly depends on $q$, and so these functions do not fall into the classes considered so far. However for every nonzero $x\in\f{q}$ we have
\[ f(x)=x^{-1}+g(x)\]
and so the equation
\[xf(x)=1+xg(x)\]
actually does not depend on $q$ anymore, a fact we will use below.
Functions of the form $f(x)=x^{q-2}+g(x)$ are in particular interesting for cryptography as another important criterion for functions used in symmetric ciphers is a high algebraic degree and functions of the studied form provide the maximal degree possible for a balanced function. Moreover, it turns out that these functions are in particular suitable for the question posed above.

The main idea to prove that a given polynomial is APN only for finitely many extensions involves the estimation of the number of points on the following surface $X'$ of the affine equation
$$
{g(x_0)+g(x_1)+g(x_2)+g(x_0+x_1+x_2)
\over (x_0+x_1)(x_2+x_1)(x_0+x_2)} x_0 x_1 x_2(x_0+x_1+x_2)=1
$$
associated to a function $f$. There are three major steps to derive these estimations.

The first, and for most of the classes of functions studied so far the most involved step, is to show that the number of points of the associated surface can be bounded by applying the Weil bound (or improvements of this bound).
The second step is the observation that the number of $\f{q}$-rational points on a certain surface associated to an APN function can be
upper bounded. The third step consists of applying the Weil bound (or improvements of this bound) to get a lower bound on the number of points
on the same surface. These bounds will involve the field size $q$. Moreover, if $q$ is large enough, the derived lower bound exceeds the upper bound for APN functions and this in turn applies that the function can only be APN up to a certain field size. This is the approach taken by all the papers that deal with this kind of questions.

For the first step, i.e. for the Weil bound to be applicable to $X'$, the surface $X'$ has to fulfill certain properties. One possibility is to show that the surface is absolutely irreducible. For most of the classes studied so far the question of absolute irreducibility was the part where one had to make further restrictions on the studied functions which where not always fully satisfactory.
However, for the functions we study the question of absolute irreducibility can be answered completely, see the Theorem below.
\begin{thm}\label{thm:absoluteirreducible}
Let $g \in \f{q}[x]$ be any non-affine polynomial. The surface $X'$ defined by the affine equation
$$
{g(x_0)+g(x_1)+g(x_2)+g(x_0+x_1+x_2)
\over (x_0+x_1)(x_2+x_1)(x_0+x_2)} x_0 x_1 x_2(x_0+x_1+x_2)=1
$$
  is absolutely irreducible.
\end{thm}
The proof of this Theorem can be found in Section \ref{sec:irr}.

The second step, i.e. the upper bound on the number of $\f{q}$-rational points on the surface $X'$ is stated in the corollary below and proven in Section \ref{sec:upperbound}.
\begin{cor}\label{fewPointsAbsIrr}
If the  polynomial mapping  $f(x)=x^{q-2}+g(x)$ ($g$ of degree $d\ge3$) is APN than the projective surface $X'$ with affine equation
$$
{g(x_0)+g(x_1)+g(x_2)+g(x_0+x_1+x_2)
\over (x_0+x_1)(x_2+x_1)(x_0+x_2)} x_0 x_1 x_2(x_0+x_1+x_2)=1
$$
has at most  $4 d q+4 q+8$ rational points.
\end{cor}
The final, and straight forward step, is a lower bound on the number of points on $X'$. One possible lower bound is given below.
\begin{cor}\label{cor:lowerboundabsoluteirreducible}
The projective surface $X'$ with affine equation
$$
{g(x_0)+g(x_1)+g(x_2)+g(x_0+x_1+x_2)
\over (x_0+x_1)(x_2+x_1)(x_0+x_2)} x_0 x_1 x_2(x_0+x_1+x_2)=1
$$
has at least $q^2+q+1-  d(d-1)  q^{3/2}-18(d+4)^4q$ rational points.
\end{cor}
\begin{pf}
As it is proven in Theorem \ref{thm:absoluteirreducible} that $X'$ is absolute irreducible we can apply a result of Lang-Weil \cite{lw} improved by Ghorpade-Lachaud (\cite[section 11]{Ghorpade99etalecohomology}), and deduce
$$|X'(\f{q})-q^2-q-1|\le d(d-1) q^{3/2}+18(d+4)^4q.$$
Hence
$$-X'(\f{q})+q^2+q+1\le d(d-1) q^{3/2}+18(d+4)^4q.$$
that is
$$X'(\f{q})\ge q^2+q+1-  d(d-1)  q^{3/2}-18(d+4)^4q.$$
\end{pf}

These results are put together in the next theorem.
\begin{thm}
\label{lawe}
Let $g$ be a polynomial from $\f{q}$ to $\f{q}$, $d$ its degree.
Then, if $d <0.45 q^{1/4} - 3.5 $ and $d\ge5$\ , $f$ is not APN.
\end{thm}

\begin{pf}
From the above it follows that if $ q^2+q+1-  d(d-1) q^{3/2}-18(d+4)^4q> 4 d q+4 q+8$, then
$X'(\f{q})>4 d q+4 q+8$, and thus $f$ is not APN.
This condition can be written as
$$ q^2-  d(d-1)q^{3/2}-(18(d+4)^4+4d+3 )q+1> 0$$
or
$$ q-  d(d-1)q^{1/2}-(18(d+4)^4+4d+3 ) +{1\over q}> 0$$
This condition is fulfilled for
$q^{1/2}>70+ 33.15 d + 4.773 d^2$.
Or for
$d <0.45 q^{1/4} - 3.5 $ and $d\ge5$.
\end{pf}

Our main result is a corollary of this.
\begin{cor} Let $g$ be any fixed non-affine polynomial function in $\f{q}[x]$. Then the
functions
\begin{eqnarray*}
f: \f{q^n} &\to& \f{q^n} \\
f(x)&=&x^{q^n-2}+g(x)
\end{eqnarray*}
are APN on at most a finite number of fields $\f{q^n}$.
\end{cor}

\subsection{Special cases}
If we fix the degree of $g$ a closer analysis of the surface $X'$ allows
to derive better bounds on the maximal field size such that the given function can be APN. If this maximal size is not too big it can be checked with the help of computers if $f$ of the given form can ever be APN by just checking up to the maximal possible extension. Following this approach to combine theoretical results and computer experiments we prove the following theorem as shown in Section \ref{sec:improvments}.
\begin{thm}\label{thm:specialcase} Let $m>3$ and $q=2^m$ be given. Furthermore, let $g$ be any non-affine polynomial of degree at most $6$ in $\f{q}[x]$. Then the
function
\begin{eqnarray*}
f: \f{q} &\to& \f{q} \\
f(x)&=&x^{q-2}+g(x)
\end{eqnarray*}
is not APN.
\end{thm}

Another special case is the case of binomials, i.e. the case where $g$ is a monomial $ax^d$. This case is particularly suitable for checking the APN property with the help of computers, as there are at most $\gcd(d,q-1)$ non equivalent functions of the form $x^{q-2}+ax^d$. This in turn implies that one can check all functions of the form at least for field sizes smaller or equal to $2^{24}$ within hours on a standard PC. For this case Theorem \ref{lawe} states that if $d$ is smaller than $30$ the function $x^{q-2}+ax^d$ can be APN only on fields up to a size of $2^{24}$. The following result therefor is again a combination of computer search on fields up to that size and Theorem \ref{lawe}.
\begin{thm}\label{thm:} Let $q=2^m$ and $d$ a nonzero integer  not a power of 2 be given. Then the
function
\begin{eqnarray*}
f: \f{q} &\to& \f{q} \\
f(x)&=&x^{q-2}+ax^d
\end{eqnarray*}
is not APN for $a\in\f q^*$ and $d\le29$.
\end{thm}

\section{Proof of Theorem \ref{thm:absoluteirreducible}: Irreducibility of $X'$}\label{sec:irr}
As stated above, one key tool is to apply the Weil bound, or variants, on the number of points on certain surfaces. In order for this bound to be applicable the surface has to fulfill certain properties. The purpose of this section is to show that the surface $X'$ defined by the affine equation
$$
{g(x_0)+g(x_1)+g(x_2)+g(x_0+x_1+x_2)
\over (x_0+x_1)(x_2+x_1)(x_0+x_2)} x_0 x_1 x_2(x_0+x_1+x_2)=1
$$
is absolutely irreducible where $g$ is any fixed polynomial.

Assume that $X'$ is not absolutely irreducible and denote
\[\phi(x_0, x_1, x_2)= \frac{g(x_0)+g(x_1)+g(x_2)+g(x_0+x_1+x_2) }{ (x_0+x_1)(x_2+x_1)(x_0+x_2)} .\]
Then one may write, with $P_i$ (resp. $Q_i$) polynomials in 2 (resp. 3) variables:
\begin{eqnarray*}
&&1+\phi(x_0, x_1, x_2){x_0 x_1 x_2(x_0+x_1+x_2)}\\
&=&(P_1(x_1, x_2)+x_0Q_1(x_0, x_1, x_2))(P_2(x_1, x_2)+x_0Q_2(x_0, x_1, x_2))\\
&=&P_1(x_1, x_2)P_2(x_1, x_2)\\
&&+x_0(Q_1(x_0, x_1, x_2)P_2(x_1, x_2)+Q_2(x_0, x_1, x_2)P_1(x_1, x_2))+\\
&&+x_0^2Q_1(x_0, x_1, x_2)Q_2(x_0, x_1, x_2)
\end{eqnarray*}
One has $P_1(x_1, x_2)P_2(x_1, x_2)=1$ so $P_1(x_1, x_2)$ is a nonzero constant, and one can suppose that $P_1(x_1, x_2)=1$ and $P_2(x_1, x_2)=1$. Consequently
\begin{eqnarray*}
&&1+\phi(x_0, x_1, x_2){x_0 x_1 x_2(x_0+x_1+x_2)}\\
&=&1+x_0(Q_1(x_0, x_1, x_2)+Q_2(x_0, x_1, x_2))+x_0^2Q_1(x_0, x_1, x_2)Q_2(x_0, x_1, x_2)\\
&=&1+x_0(Q_1(0, x_1, x_2)+Q_2(0, x_1, x_2))+\\
&&x_0^2\pt{{Q_1(0, x_1, x_2)+Q_2(0, x_1, x_2)+Q_1(x_0, x_1, x_2)+Q_2(x_0, x_1, x_2)\over x_0}}+\\
&&x_0^2Q_1(x_0, x_1, x_2)Q_2(x_0, x_1, x_2)
\end{eqnarray*}
and
\begin{eqnarray*}
&&1+\phi(x_0, x_1, x_2){x_0 x_1 x_2(x_0+x_1+x_2)}\\
&=&1+\phi(0, x_1, x_2){x_0 x_1 x_2(x_1+x_2)}+\phi(0, x_1, x_2){x_0^2 x_1 x_2}+\\
&&{x_0^2 x_1 x_2(x_0+x_1+x_2)}\pt{{\phi(x_0, x_1, x_2)+\phi(0, x_1, x_2)\over x_0}}
\end{eqnarray*}
where the fractions with denominators $x_0$ are actually polynomials.
Hence
$$Q_1(0, x_1, x_2)+Q_2(0, x_1, x_2)=\phi(0, x_1, x_2){x_1 x_2(x_1+x_2)}$$
and
\begin{eqnarray*}
&&{Q_1(0, x_1, x_2)+Q_2(0, x_1, x_2)+Q_1(x_0, x_1, x_2)+Q_2(x_0, x_1, x_2)\over x_0}+\\
&&Q_1(x_0, x_1, x_2)Q_2(x_0, x_1, x_2)\\
&=&{x_1 x_2(x_0+x_1+x_2)}{\phi(x_0, x_1, x_2)+\phi(0, x_1, x_2)\over x_0}+\phi(0, x_1, x_2){ x_1 x_2}
\end{eqnarray*}
Let $d$ be the degree of $g$.
Remark that the degree of $ \phi(0, x_1, x_2)$ is equal to $d-3$ (see Lemma \ref{degree}).
Hence
$$d=\deg (\phi(0, x_1, x_2){x_1 x_2(x_1+x_2)})\le \sup(\deg Q_1,Q_2)$$
and
$\deg Q_1+\deg Q_2\le 3+d-3-1$.
If $\deg Q_1\ge\deg Q_2$, one has
$$d\le \sup(\deg Q_1,Q_2)\le\deg Q_1\le \deg Q_1+\deg Q_2\le d-1$$
We obtain a contradiction to the assumption that $X'$ is reducible.

\section{Proof of Corollary \ref{fewPointsAbsIrr}: The Upper Bound}\label{sec:upperbound}
The purpose of this section is to give a proof of Corollary \ref{fewPointsAbsIrr}, i.e. to show that the number of rational
points on the surface $X'$ can be upper bounded if $f$ is APN. The main tool is the following Lemma (see for example \cite{rodier2006bdp}).

Let $f$ be a a polynomial mapping from $\f{q}$ to itself which has no terms of degree a power of 2.
\begin{prop}\label{fewPoints}
\label{apn}
The function
$ f:\f{q}\to \f{q}$
is APN if and only if the surface
$ f(x_0)+ f(x_1)+ f(x_2)+ f(x_0+x_1+x_2)=0$
has all of its rational points contained in the surface
$(x_0+x_1)(x_2+x_1)(x_0+x_2)=0$.
\end{prop}

Before we prove Corollary \ref{fewPointsAbsIrr}, remark that
the polynomial $f(x_0)+f(x_1)+f(x_2)+f(x_0+x_1+x_2)$ is divisible by
$(x_0+x_1)(x_2+x_1)(x_0+x_2)$, therefore the quotient
$${f(x_0)+f(x_1)+f(x_2)+f(x_0+x_1+x_2)\over (x_0+x_1)(x_2+x_1)(x_0+x_2)}$$
defines a polynomial which is the affine equation of a surface $X$ if the polynomial is not constant,that is if $f$ is not a $q$-affine  polynomial or a polynomial of degree 3.

We will make use of the following lemmata.
\begin{lem}
If $d$ is not a power of 2 and at least 3, and an integer $c$ at least 2  then the polynomial
$\sum_{i=1}^c x_i^d +(\sum_{i=1}^c x_i)^d$ is non zero, hence of degree $d$.
\end{lem}
\begin{pf}
Write $d=2^ab$ with $b$ odd.
The polynomial can be written
$(\sum_{i=1}^c x_i^b +(\sum_{i=1}^c x_i)^b)^{2^a}$.
The inner polynomial contains a monomial (say $bx_1x_2^{b-1}$) of degree $b$.
Hence the conclusion.
\end{pf}

\begin{lem}\label{degree}
If $\deg g=d$ is not a power of 2, then
$$\phi(x_0,x_1,x_2)={g(x_0)+g(x_1)+g(x_2)+g(x_0+x_1+x_2)\over (x_0+x_1)(x_2+x_1)(x_0+x_2)}$$
is a polynomial  of degree $d-3$.
\end{lem}
\begin{pf}
Denote by $\phi_{d}$ the term of highest degree of $\phi$ . As it is nonzero and as $\phi_{d}$ is a rational homogeneous fraction, $\phi_{d}$ is of degree $d-3$.
\end{pf}

Now we can prove Corollary \ref{fewPointsAbsIrr}:
Let the  polynomial mapping $f(x)=x^{q-2}+g(x)$ ($g$ of degree $d\ge3$) be APN and $X'$ be the surface with affine equation
$$
{g(x_0)+g(x_1)+g(x_2)+g(x_0+x_1+x_2)
\over (x_0+x_1)(x_2+x_1)(x_0+x_2)} x_0 x_1 x_2(x_0+x_1+x_2)=1
$$
Due to Theorem \ref{thm:absoluteirreducible} $X'$ is absolutely irreducible. We have to show that the corresponding projective surface has at most  $4 d q+4 q+8$ rational points.

If the surface $X'$ contained the plane $x_0+x_1=0$, it would contain also the planes $x_2+x_1=0$ and $x_0+x_2=0$ by symmetry, which is impossible as the surface $X'$ is irreducible.
So its intersection with the plane $x_0+x_1=0$ is a curve of degree $d+1$. This curve has at most $(d+1)q+1$ rational points from Serre \cite{se}. The same argument works for the plane at l'infinity.

If $f$ is APN, the affine surface $X$ has no other rational points than those of the surface $(x_0+x_1)(x_2+x_1)(x_0+x_2)=0$, which is union of a plane $x_0+x_1=0$  and of its symmetrical plane.

The set $\cj X'(\f q)$ decomposes as follows:
$$\cj X'(\f q)=X'_{x_0}\cup X'_{x_1}\cup X'_{x_2}\cup X'_{x_0+x_1+x_2}\cup X^*_{\hbox {\tiny aff}}
\cup X^*_\infty$$
where $X'_{a}=\cj X'(\f q)\cap(a=0)$, and $X^*_{\hbox {\tiny aff}}$ is the affine complement.
The equation of  surface $X'$ may be written as follows
$$\displaylines{
{x_0^{-1}+g(x_0)+x_1^{-1}+g(x_1)+x_2^{-1}+g(x_2)+(x_0+x_1+x_2)^{-1}+g(x_0+x_1+x_2)
\over (x_0+x_1)(x_2+x_1)(x_0+x_2)}\hfill\cr\hfill
\times x_0 x_1 x_2(x_0+x_1+x_2)=0.
}$$
It means that for $x_0 x_1 x_2(x_0+x_1+x_2)\ne0$, the element of the set $X^*_{\hbox {\tiny aff}}$ fulfill
$$\displaylines{
{f(x_0)+f(x_1)+f(x_2)+f(x_0+x_1+x_2)
\over (x_0+x_1)(x_2+x_1)(x_0+x_2)}
\times x_0 x_1 x_2(x_0+x_1+x_2)=0
}$$
which proves that $X^*_{\hbox {\tiny aff}}$ is contained in $X(\f q)$ hence, as $f$ is APN, in the union of the three planes $(x_0+x_1)(x_2+x_1)(x_0+x_2)=0$.
Therefore the number of points in  $X^*_{\hbox {\tiny aff}}$ is bounded by $3((d+1)q+1)$.

The intersection of the surface $\cj X'$ with the plane $x_0=0$ is the line $x_0=0$ in the plane at infinity, which has  $q+1$ rational points.

The equation of the intersection of the surface $\cj X'$ with the plane at infinity is
$$
{g(x_0)+g(x_1)+g(x_2)+g(x_0+x_1+x_2)
\over (x_0+x_1)(x_2+x_1)(x_0+x_2)} x_0 x_1 x_2(x_0+x_1+x_2)=0
$$
It contains the lines $x_0 x_1 x_2(x_0+x_1+x_2)=0$, for which we have already taken in consideration the points, and the curve
$$
{g(x_0)+g(x_1)+g(x_2)+g(x_0+x_1+x_2)
\over (x_0+x_1)(x_2+x_1)(x_0+x_2)}=0
$$
which has at most $(d-3)q+1$ rational points.

So
\begin{eqnarray*}
\#\cj X'(\f q)&=&\#X'_{x_0} +  \#X'_{x_1} +  \#X'_{x_2} +  \#X'_{x_0+x_1+x_2} +  X^*_{\hbox {\tiny aff}}
 +  \#X^*_\infty\\
&=&4(q+1)+  3((d+1)q+1)+(d-3)q+1\\
&=&4 d q+4 q+8
 \end{eqnarray*}
which proves the result.

\section{Improvements in specific cases}\label{sec:improvments}
Under certain conditions, we can obtain a better bound on the dimension.

\begin{thm}
\label{lisse}
Let $g$ a polynomial mapping  from $\f{2^m}$ to itself, $d$ its degree.  Let us suppose that the surface $X'$ defined by
$$
{g(x_0)+g(x_1)+g(x_2)+g(x_0+x_1+x_2)
\over (x_0+x_1)(x_2+x_1)(x_0+x_2)} x_0 x_1 x_2(x_0+x_1+x_2)=1
$$
of degree $d'=d+1$ has only isolated singular points.
Then if  $d\ge3$  and $d <q^{1/4} $\ , $f$ is not APN.
\end{thm}

\begin{pf}

From an improvement of a result  of Deligne \cite{Deligne} by Ghorpade-Lachaud (\cite{Ghorpade99etalecohomology}, corollaire 7.2), we deduce
\begin{eqnarray*}
 |X'(\f{q})-q^2-q-1|
 &\le& b'_{1}(2, d') q^{3/2} + (b_2(3, d') + 1) q\\
 &\le& (d'-1)(d'-2) q^{3/2} + (d'^3-4d'^2+6d' - 1) q\\
 &\le& d(d-1) q^{3/2} + (2 + d - d^2 + d^3) q
\end{eqnarray*}
where $b'_{1}$ and $b_2(3, d')$ are Betti numbers (see \cite{Ghorpade99etalecohomology}).
         Hence
         $$X'(\f{q})\ge q^2+q+1- d(d-1) q^{3/2} - (2 + d - d^2 + d^3) q$$

Therefore if
$$q^2+q+1- d(d-1) q^{3/2} - (2 + d - d^2 + d^3) q> 4dq+4q+8,$$ then
$X'(\f{q})>4dq+4q+8$, and $f$ is not APN due to Corollary \ref{fewPointsAbsIrr}.
This condition can be rewritten as
$$q- d(d-1) q^{1/2} - (5 +5 d - d^2 + d^3) -7/q> 0.$$
It is fulfilled if $q>d^4$ as soon as $d\ge3$.
\end{pf}

\subsection{Polynomials $g$ of small degree}

As stated in the introduction the APN property is invariant under the so called CCZ-equivalence. As adding affine functions is a special case of CCZ-equivalence it is clear that given two functions $f$ and $g$ that differ by an affine function $f$ is APN if and only if $g$ is APN. Moreover, multiplying an APN polynomial function by a constant or replacing $x$ by any non-constant linear polynomial yields again an APN polynomial.
These well known observations are summarized in the proposition below.
\begin{prop}
\label{p2}
A polynomial function  $f$ is APN if and only if the polynomial $f_0$ obtained by removing all monomials of degree a power of 2 and by removing the constant term  is APN.
Moreover, if $f$ is APN then for any nonzero $a,c\in \f{2^m}$ and any element $b\in \f{2^m}$ the polynomial function
\[ cf(ax+b) \]
is APN.
\end{prop}
This proposition will be used in the reminder of the paper to simplify the polynomials we have to study.

\paragraph{Polynomials of degree 3}
We first focus on polynomials of degree $3$. Here, the general form of $f$ is
\[ f(x)=x^{q-2}+a_3x^3+a_2x^2+a_1x+a_0\]
where $a_3 \ne 0$ which is clearly equivalent to
\[ f(x)=x^{q-2}+a_3x^3.\]
Moreover, replacing $x$ by $a_3^{-1/4}x$ and multiplying across by $a_3^{1/4}$ we see that $f$ is in any case equivalent to
\[f(x)=x^{q-2}+x^3.\]
In this case the affine equation for $X'$ becomes

\begin{eqnarray*}
&&1+{g(x_0)+g(x_1)+g(x_2)+g(x_0+x_1+x_2)
\over (x_0+x_1)(x_2+x_1)(x_0+x_2)} x_0 x_1 x_2(x_0+x_1+x_2) \\
&=&1+{\pt{(x_0^3+x_1^3+x_2^3+(x_0+x_1+x_2)^3)x_0 x_1 x_2(x_0+x_1+x_2)}
\over (x_0+x_1)(x_2+x_1)(x_0+x_2)}\\
&=&1+{x_0 x_1 x_2(x_0+x_1+x_2)}\\
\end{eqnarray*}
The search of singular points on the surface
$z^4+{x_0 x_1 x_2(x_0+x_1+x_2)}=0$
gives a finite number of points.

\paragraph{Polynomials of degree $5$}
Next, we study polynomials of degree $5$. The general form of $f$ (up to adding affine equivalence) is given by
\[ f(x)=x^{q-2}+a_5x^5+a_3x^3.\]
Furthermore we can assume without loss of generality that $x_3 \in \f{2}$.
We have to study the surface
\begin{eqnarray*}
    X'&=1+&\left(a_3\frac{x_0^3+x_1^3+x_2^3+(x_0+x_1+x_2)^3}{(x_0+x_1)(x_0+x_2)(x_1+x_2)} \right.\\
    &&\left.+a_5\frac{x_0^9+x_1^3+x_2^9+(x_0+x_1+x_2)^9}{(x_0+x_1)(x_0+x_2)(x_1+x_2)} \right)x_0x_1x_2(x_0+x_1+x_2)
\end{eqnarray*}
and show that there are only a finite number of singular points. The lengthy -- but straight forward -- computation for showing this can be found in Appendix \ref{Ap:Deg5}.

\paragraph{Polynomials of degree $6$}
Next, we study polynomials of degree $6$. The general form of $f$ (up to adding affine equivalence) is given by
\[ f(x)=x^{q-2}+a_6x^6+a_5x^5+a_3x^3.\]
Furthermore we can assume without loss of generality that $x_3 \in \f{2}$.

In this case we have to study the surface
\begin{eqnarray*}
    X'&=1+&\left(a_3\frac{x_0^3+x_1^3+x_2^3+(x_0+x_1+x_2)^3}{(x_0+x_1)(x_0+x_2)(x_1+x_2)} \right.\\
    &&+a_6\frac{x_0^6+x_1^3+x_2^6+(x_0+x_1+x_2)^6}{(x_0+x_1)(x_0+x_2)(x_1+x_2)}\\
    &&\left.+a_5\frac{x_0^9+x_1^3+x_2^9+(x_0+x_1+x_2)^9}{(x_0+x_1)(x_0+x_2)(x_1+x_2)} \right)x_0x_1x_2(x_0+x_1+x_2)
\end{eqnarray*}
and show that there are only a finite number of singular points. We refer to Appendix \ref{Ap:Deg6} for the proof.
\paragraph{Conclusion}
As we have seen above for any non-affine polynomial of degree less than 6 the corresponding surface contains only isolated singularities. Therefore, Proposition \ref{p2} applies. Thus if $f$ is APN it holds that $q=2^m \le 6^4$ which implies $m \le 10$.
It can easily be checked that functions form $x^{q-2}+a_6x^6+a_5x^5+a_3x^3$ are APN only if $m \le 3$. Note that, for $m\le 3$ the APN functions are quadratic and moreover for $m \le 3$ all APN functions are CCZ equivalent to $x^3$. These considerations finally prove Theorem \ref{thm:specialcase} stated above.
\bibliographystyle{plain}

\begin{thebibliography}{10}

\bibitem{bracken2008fmq}
C.~Bracken, E.~Byrne, N.~Markin, and G.~McGuire.
\newblock A few more quadratic apn functions.
\newblock {\em Arxiv preprint arXiv:0804.4799}, 2008.

\bibitem{bracken2007nfq}
C.~Bracken, E.~Byrne, N.~Markin, and G.~McGuire.
\newblock New families of quadratic almost perfect nonlinear trinomials and
  multinomials.
\newblock {\em Finite Fields and Their Applications}, 14(3):703--714, 2008.

\bibitem{budaghyan2008cqa}
L.~Budaghyan and C.~Carlet.
\newblock Classes of quadratic apn trinomials and hexanomials and related
  structures.
\newblock {\em Information Theory, IEEE Transactions on}, 54(5):2354--2357,
  2008.

\bibitem{TraceAPN}
L.~Budaghyan, C.~Carlet, and G.Leander.
\newblock Constructing new apn functions from known ones.
\newblock Finite Fields and Applications, to appear.

\bibitem{budaghyan2008tcq}
L.~Budaghyan, C.~Carlet, and G.~Leander.
\newblock Two classes of quadratic apn binomials inequivalent to power
  functions.
\newblock {\em Information Theory, IEEE Transactions on}, 54(9):4218--4229,
  2008.

\bibitem{DBLP:journals/dcc/CarletCZ98}
Claude Carlet, Pascale Charpin, and Victor Zinoviev.
\newblock Codes, bent functions and permutations suitable for {DES}-like
  cryptosystems.
\newblock {\em Des. Codes Cryptography}, 15(2):125--156, 1998.

\bibitem{Deligne}
Pierre Deligne.
\newblock {La conjecture de Weil : I.}
\newblock {\em Publications Mathema\-tiques of l'IHES}, 43:273--307, 1974.

\bibitem{DBLP:journals/IEEE/Pot06}
Y.~Edel, G.~Kyureghyan, and A.~Pott.
\newblock A new apn function which is not equivalent to a power mapping.
\newblock {\em IEEE Transactions on Information Theory}, 52(2):744--747, 2006.

\bibitem{Ghorpade99etalecohomology}
Sudhir~R. Ghorpade and Gilles Lachaud.
\newblock Etale cohomology, {Lefschetz} theorems and the number of points of
  singular varieties over finite fields.
\newblock {\em Moscow Mathematical Journal}, 2:589--631, 2002.

\bibitem{Jedlicka}
David Jedlicka.
\newblock Apn monomials over {${\rm GF}(2\sp n)$} for infinitely many $n$.
\newblock {\em Finite Fields and Their Applications}, 13(4):1006--1028, 2007.

\bibitem{lw}
Serge Lang and Andre Weil.
\newblock Number of points of varieties in finite fields.
\newblock {\em American Journal of Mathematics}, 76(4):819--827, 1954.

\bibitem{rodier2006bdp}
F.~Rodier.
\newblock {Borne sur le degr\'e des polyn\^omes presque parfaitement
  non-lin\'eaires}.
\newblock {\em Arxiv preprint math.AG/0605232, to be published with the
  proceedings of the conference AGCT-11}, 2006.

\bibitem{rodierBFCA08}
F.~Rodier.
\newblock Bounds on the degrees of apn polynomials.
\newblock {\em to be published with the proceedings of the workshop BFCA08,
  Copenhagen, 2008}, 2006.

\bibitem{se}
J.~P. Serre.
\newblock Lettre \`a {M. Tsfasman}.
\newblock {\em Asterisque}, 198-199-200:351--353, 1991.

\bibitem{Voloch}
Felipe Voloch.
\newblock Symmetric cryptography and algebraic curves.
\newblock Proceedings of the First SAGA Conference, Papeete, France, 2007.

\end{thebibliography}

\begin{appendix}

\section{Singular points for $g$ of degree $5$}\label{Ap:Deg5}
We have to study the surface
\begin{eqnarray*}
    X'&=1+&\left(a_3\frac{x_0^3+x_1^3+x_2^3+(x_0+x_1+x_2)^3}{(x_0+x_1)(x_0+x_2)(x_1+x_2)} \right.\\
    &&\left.+a_5\frac{x_0^9+x_1^3+x_2^9+(x_0+x_1+x_2)^9}{(x_0+x_1)(x_0+x_2)(x_1+x_2)} \right)x_0x_1x_2(x_0+x_1+x_2)
\end{eqnarray*}
and show that there are only a finite number of singular points.
For this we compute the derivatives of the projective version of $X'$
\begin{eqnarray*}
X'&=&  a_5x_2x_1x_0(x_0 + x_1 + x_2)(x_0^2 + x_0x_1 + x_0x_2 + x_1^2 + x_1x_2 + x_2^2)\\
&& +  a_3z^2x_2x_1x_0(x_0 + x_1 + x_2)\\
&&+z^6
\end{eqnarray*}
The derivatives of the projective version of $X'$ are as follows
\begin{eqnarray*}
    \frac{\partial X'}{\partial x_0}&=&x_2x_1(x_1+x_2)P(x_1,x_2,z) \\
    \frac{\partial X'}{\partial x_1}&=&x_0x_2(x_0+x_2)P(x_0,x_2,z) \\
    \frac{\partial X'}{\partial x_2}&=&x_0x_1(x_0+x_1)P(x_0,x_1,z) \\
    \frac{\partial X'}{\partial z}&=& 0
\end{eqnarray*}
where
\[P(x,y,z)=a_5(x^2 +xy + y^2) + a_3z^2 \]
To study the singular points of these equations, we make some case distinction.
\paragraph{Case $x_0=0$:} In this case $X'$ simplifies to $z=0$ and we have
\[ \frac{\partial X'}{\partial x_0}(0,x_1,x_2,0)=a_5x_2x_1(x_1+x_2)(x_1^2 +x_1x_2 + x_2^2) \]
which, up to equivalence, implies a finite number of singularities .
\paragraph{Case $x_1=0$ or $x_2=0$:}
Due to symmetries the cases can be handled exactly like the first case.
\paragraph{Case $x_1=x_2$:}
Here we are left with the following system of equations
\begin{eqnarray*}
a_5x_2^2x_0^2 (x_0 + x_2 )^2 +  a_3z^2x_2^2x_0^2+z^6&=&0 \\
x_0x_2(x_0+x_2)(a_5(x_0^2 +x_0x_2 + x_2^2) + a_3z^2) &=&0
\end{eqnarray*}
Now if $x_0=x_2$ the first of this equation becomes
\begin{eqnarray}\label{usedAgainForDeg6}
 (a_3x_0^4+z^4)z^2&=&0
\end{eqnarray}
If $a_3=0$ then $z=0$ and there are, up to equivalence, at most two points $(1,1,1,0)$ and $(0,0,0,0)$. For $a_3 \ne 0$ we can assume $a_3=1$, see above. Now, if $z\ne 0$ then $x_0=z$ and, up to equivalence, there is only one point $(1,1,1,1)$.

The case $z=0$ and $x_0=x_1=x_2$ correspond, up to equivalence to at most two points  $(1,1,1,0)$ and $(0,0,0,0)$.

For $x_0=0$ or $x_2=0$ see the first case above. Thus we are left with the case
\begin{eqnarray*}
a_5(x_0^2 +x_0x_2 + x_2^2) + a_3z^2 &=&0 \\
a_5x_2^2x_0^2(x_0 + x_2 )^2 + a_3x_0^2x_2^2z^2 + z^6&=&0
\end{eqnarray*}
Computing the resultants of both equations with respect to $x_0$ and $x_2$ we get
\[ a_5^3Q(x_0)=0 \mbox{ and } a_5^3Q(x_2)=0 \]
where
\begin{eqnarray*}
    Q(x)&=& a_5^3x^{12} + a_3a_5^2x^{10}z^2 + a_3^2a_5x^8z^4 + a_3^3x^6z^6 +
     a_3a_5x^4z^8 + a_5z^{12}.
\end{eqnarray*}
The fact that $Q$ is a non-zero polynomial (as $a_5 \ne 0$) implies only a finite number of solutions (up to equivalence).
\paragraph{Case $x_0=x_2$ or $x_0=x_1$:}
Due to symmetries the cases can be handled exactly like the case above.
\paragraph{None of the above cases:}
Here we get the following equations.
\begin{eqnarray*}
X'&=&0 \\
P(x_1,x_2,z)&=&0 \\
P(x_0,x_2,z)&=&0 \\
P(x_0,x_1,z)&=&0
\end{eqnarray*}
In particular the singular points have to fulfil
\begin{eqnarray*}
P(x_1,x_2,z)+P(x_0,x_2,z)&=& a_5(x_0+x_1)(x_0+x_1+x_2) =0  \\
P(x_1,x_2,z)+P(x_0,x_2,z)&=& a_5(x_0+x_2)(x_0+x_1+x_2) =0 \\
P(x_1,x_2,z)+P(x_0,x_2,z)&=& a_5(x_1+x_2)(x_0+x_1+x_2) =0.
\end{eqnarray*}
The cases $x_0=x_1$ or $x_0=x_2$ or $x_1=x_2$ have been studied above. The case $x_0+x_1+x_2=0$ implies $z=0$ (as $X'=0$) and
then
\[P(x_0,x_2,0)=a_5(x_0^2 +x_0x_1 + x_1^2) =0\]
which leaves only a finite number of singularities (up to equivalence).

\section{Singular points for $g$ of degree $6$}\label{Ap:Deg6}
We have to study the surface
\begin{eqnarray*}
    X'&=1+&\left(a_3\frac{x_0^3+x_1^3+x_2^3+(x_0+x_1+x_2)^3}{(x_0+x_1)(x_0+x_2)(x_1+x_2)} \right.\\
    &&+a_6\frac{x_0^6+x_1^3+x_2^6+(x_0+x_1+x_2)^6}{(x_0+x_1)(x_0+x_2)(x_1+x_2)}\\
    &&\left.+a_5\frac{x_0^9+x_1^3+x_2^9+(x_0+x_1+x_2)^9}{(x_0+x_1)(x_0+x_2)(x_1+x_2)} \right)x_0x_1x_2(x_0+x_1+x_2)
\end{eqnarray*}
and show that there are only a finite number of singular points. The computations and case distinctions are very similar to the case where $g$ is of degree $5$.

We compute the derivatives of the projective version of $X'$
\begin{eqnarray*}
X'&=&  a_6x_2x_1x_0(x_1 + x_2)(x_0 + x_2)(x_0 + x_1)(x_0 + x_1 + x_2)\\
&& +  a_5zx_2x_1x_0(x_0 + x_1 + x_2)(x_0^2 + x_0x_1 + x_0x_2 + x_1^2 + x_1x_2 + x_2^2)\\
&& +  a_3z^3x_2x_1x_0(x_0 + x_1 + x_2)\\
&&+z^7
\end{eqnarray*}

with respect to $x_0,x_1,x_2$ and $z$. We get
\begin{eqnarray*}
    \frac{\partial X'}{\partial x_0}&=&x_2x_1(x_1+x_2)P(x_1,x_2,z) \\
    \frac{\partial X'}{\partial x_1}&=&x_0x_2(x_0+x_2)P(x_0,x_2,z) \\
    \frac{\partial X'}{\partial x_2}&=&x_0x_1(x_0+x_1)P(x_0,x_1,z) \\
    \frac{\partial X'}{\partial z}&=& a_5x_2x_1x_0(x_0 + x_1 + x_2)(x_0^2 + x_0x_1 + x_0x_2 + x_1^2 + x_1x_2 + x_2^2)\\
&& +  a_3z^2x_2x_1x_0(x_0 + x_1 + x_2)+z^6
\end{eqnarray*}
where
\[P(x,y,z)=a_6xy(x+y) + a_5z(x^2 +xy + y^2) + a_3z^3 \]
To study the singular points of these equations, we make some case distinction.
\paragraph{Case $x_0=0$:}
Then $\frac{\partial X'}{\partial z}(0,x_1,x_2)=0$ implies $z=0$ and this simplifies to
\begin{eqnarray*}
\frac{\partial X'}{\partial x_0}(0,x_1,x_2,0)&=&a_6x_2^2x_1^2(x_1+x_2)^2 \\
&=&0
\end{eqnarray*}
which for $a_6\ne 0$ implies, up to equivalence, a finite number of singularities.
\paragraph{Case $x_1=0$ or $x_2=0$:}
Due to symmetries the cases can be handled exactly like the first case.
\paragraph{Case $x_1=x_2$:}
In this case we are left with only two non-zero equations, namely
\begin{eqnarray*}
\frac{\partial X'}{\partial x_1}(x_0,x_2,x_2,z)&=&  x_0x_2(x_0+x_2)P(x_0,x_2,z) \\
\frac{\partial X'}{\partial z}(x_0,x_2,x_2,z)&=&  a_5x_2^2x_0^2(x_0 + x_2 )^2 + a_3x_0^2x_2^2z^2 + z^6
\end{eqnarray*}
Now, if $x_0=x_2$ then $x_0=x_1=x_2$  and we get
\[ \frac{\partial X'}{\partial z}(x_0,x_0,x_0,z)=(a_3x_0^4+z^4)z^2. \]
which has been studied already in the case of a degree 5 polynomial (see Equation (\ref{usedAgainForDeg6}))

For $x_0=0$ or $x_2=0$ see the first case above. Thus we are left with the case
\begin{eqnarray*}
P(x_0,x_2,z)&=&0 \\
a_5x_2^2x_0^2(x_0 + x_2 )^2 + a_3x_0^2x_2^2z^2 + z^6&=&0
\end{eqnarray*}
For this we again distinguish two cases: $z=0$ and $z\ne 0$. For $z=0$ we get
\[ P(x_0,x_2,0)= a_6x_2x_0(x_0+x_2) .\]
Thus we have either $x_0=0$ or $x_2=0$ which have been handled above, or $x_0=x_2$ which gives again at most two points up to equivalence, see above.

For $z\ne 0$ we can restrict, up to equivalence, to the case $z=1$. In this case we get
\begin{eqnarray*}
    a_6x_0x_2(x_0+x_2) + a_5(x_0^2 +x_0x_2 + x_2^2) + a_3 &=&0\\
    a_5x_2^2x_0^2(x_0 + x_2 )^2 + a_3x_0^2x_2^2 + 1 &=&0
\end{eqnarray*}
Computing the resultant of these two equations with respect to $x_0$ and $x_2$ we get
\[ Q(x_0)=0 \mbox{ and } Q(x_2)=0 \]
where
\begin{eqnarray*}
    Q(x)&=&(a_3a_5^3a_6^2 + a_5^6)x^{12} + (a_3^2a_5^2a_6^2 + a_3a_5^5)x^{10} +
   (a_3^3a_5a_6^2 + a_3^2a_5^4 + a_3a_6^4)x^8 \\
   &&+ (a_3^4a_6^2 + a_3^3a_5^3)x^6
   + (a_3a_5^4 + a_6^4)x^4 + a_5^4
\end{eqnarray*}
For $a_5\ne 0$ $Q$ is a non-zero polynomial as its constant term is non-zero. For $a_5=0$ $Q$ is non-zero as the degree $4$ term is non-zero. Therefore, in any case we get at most finitely many points.

\paragraph{Case $x_0=x_2$ or $x_0=x_1$:}
Due to symmetries the cases can be handled exactly like the case above.
\paragraph{None of the above cases:}
In this case we have to study
\begin{eqnarray*}
    P(x_1,x_2,z)&=&0 \\
    P(x_0,x_2,z)&=&0 \\
    P(x_0,x_1,z)&=&0 \\
    \frac{\partial X'}{\partial z}&=& 0
\end{eqnarray*}
This implies that $x_0,x_1,x_2,z$ are solutions to the following set of equations:
\begin{eqnarray*}
P(x_1,x_2,z)+P(x_0,x_2,z)&=& (a_6x_2 + a_5z)(x_0+x_1)(x_0+x_1+x_2) =0  \\
P(x_1,x_2,z)+P(x_0,x_2,z)&=& (a_6x_1 + a_5z)(x_0+x_2)(x_0+x_1+x_2) =0 \\
P(x_1,x_2,z)+P(x_0,x_2,z)&=& (a_6x_0 + a_5z)(x_1+x_2)(x_0+x_1+x_2) =0
\end{eqnarray*}
The cases $x_0=x_1$ or $x_0=x_2$ or $x_1=x_2$ have been studied above. The case $x_0+x_1+x_2=0$ implies $z=0$ (as $X'=0$) and
as seen above
\[P(x_0,x_2,0)=a_6x_2x_0(x_0+x_2) ,\]
and we are back to cases studied before. Thus the only case left is
\begin{eqnarray*}
(a_6x_2 + a_5z)&=&0  \\
(a_6x_1 + a_5z)&=&0  \\
(a_6x_0 + a_5z)&=&0
\end{eqnarray*}
which implies $x_0=x_1=x_2$ and, up to equivalence, at most two points.

\end{appendix}
\end{document}